\documentclass[12pt]{amsart}
\usepackage{graphicx,amsfonts}

\newcommand{\Q}{\mathbb Q} 
\newcommand{\R}{\mathbb R} 
\newcommand{\C}{\mathbb C} 

\newcommand{\A}{\mathcal{A}} 
\newcommand{\B}{\mathcal{B}} 

\newcommand{\eps}{\varepsilon}

\newcommand{\cdWW}{\raisebox{-2mm}{\includegraphics[height=6mm]{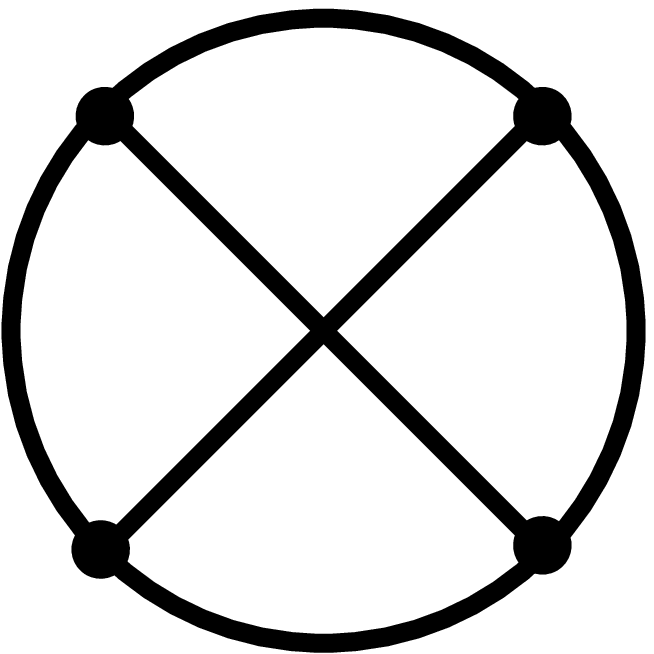}}}

\title{The Kontsevich integral}
\author{S.~Chmutov, S.~Duzhin}
\thanks{
Date: January 3, 2005.\\
\indent
The second author was partially supported by grant
RFBR NSh--1972.2003.1}

\begin{document}
\maketitle

\begin{abstract}
This is an overview article on the Kontsevich integral written for the 
\textit{Encyclopedia of Mathematical Physics} to be published by Elsevier.
\end{abstract}

\section{Introduction}
\label{intr}

The Kontsevich integral was invented by M.~Kontsevich \cite{Kon}
as a tool to prove the fundamental theorem of the theory of
finite type (Vassiliev) invariants (see \cite{BN1,BNe}).
It provides an invariant exactly as strong as the totality
of all Vassiliev knot invariants.

The Kontsevich integral is defined for oriented tangles (either framed or
unframed) in $\R^3$, therefore it is also defined in the particular cases
of knots, links and braids.

$$
  \includegraphics[width=12cm]{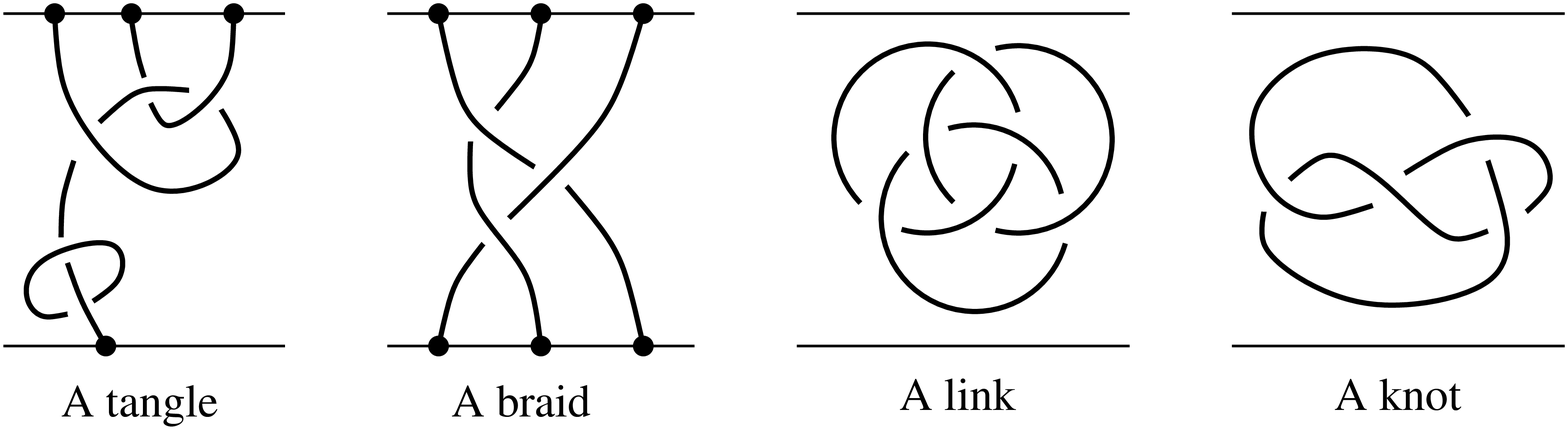}
$$

As a starter, we give two examples where simple versions of the Kontsevich
integral have a straightforward geometrical meaning. In these examples, as
well as in the general construction of the Kontsevich integral, we represent
3-space $\R^3$ as the product of a real line $\R$ with coordinate $t$ and a
complex plane $\C$ with complex coordinate $z$.
\smallskip

\noindent
\begin{minipage}{4.5in}
\textbf{Example 1.} The number of twists in a braid with two strings $z_1(t)$ 
and $z_2(t)$ placed in the slice $0\le t\le 1$ is equal to
$$\frac{1}{2\pi i}\int_0^1\frac{dz_1-dz_2}{z_1-z_2}.$$
\end{minipage}
\quad
\raisebox{-12mm}{\includegraphics[height=25mm]{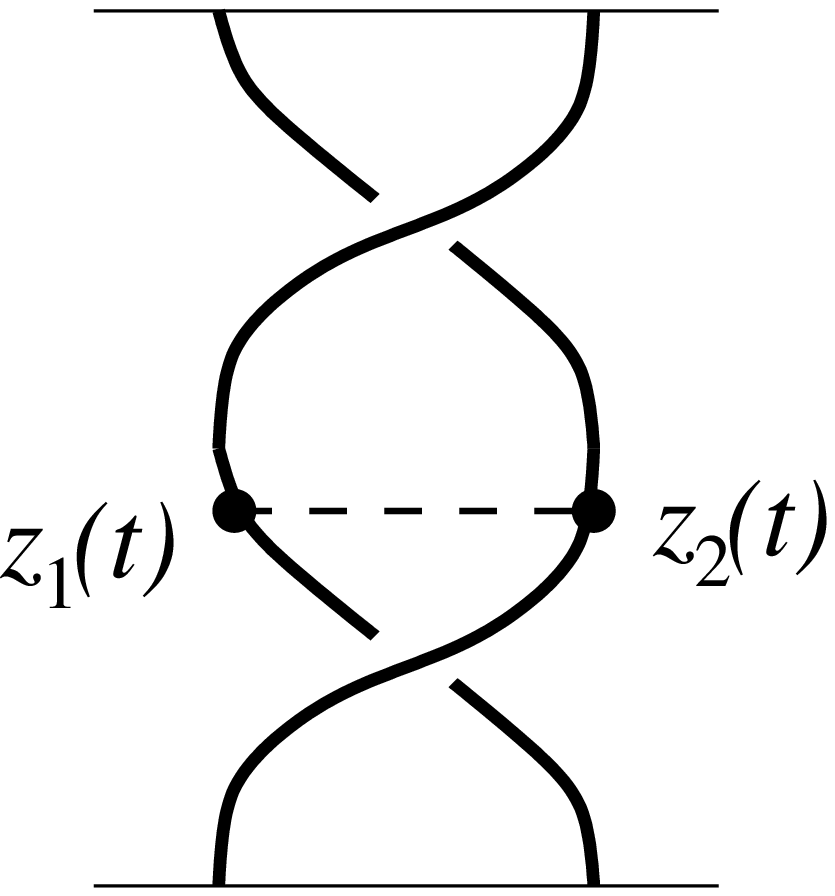}}

\medskip

\noindent
\begin{minipage}{4.5in}
\textbf{Example 2.} The linking number of two spatial curves
$K$ and $K'$ can be computed as
$$  
  lk(K,K')=\frac{1}{2\pi i}
          \int_{m<t<M}\sum_j \eps_j
          \frac{d(z_j(t)-z'_j(t))}{z_j(t)-z'_j(t)},
$$
where $m$ and $M$ are the minimum and the maximum values of $t$ on the
link $K\cup K'$, $j$ is the index that enumerates all possible choices
of a pair of strands of the link as functions $z_j(t)$, $z'_j(t)$ corresponding
to $K$ and $K'$, respectively, and $\eps_j=\pm1$ according to the parity
of the number of chosen strands that are oriented downwards.
\end{minipage}
\quad
\raisebox{-17mm}{\includegraphics[height=35mm]{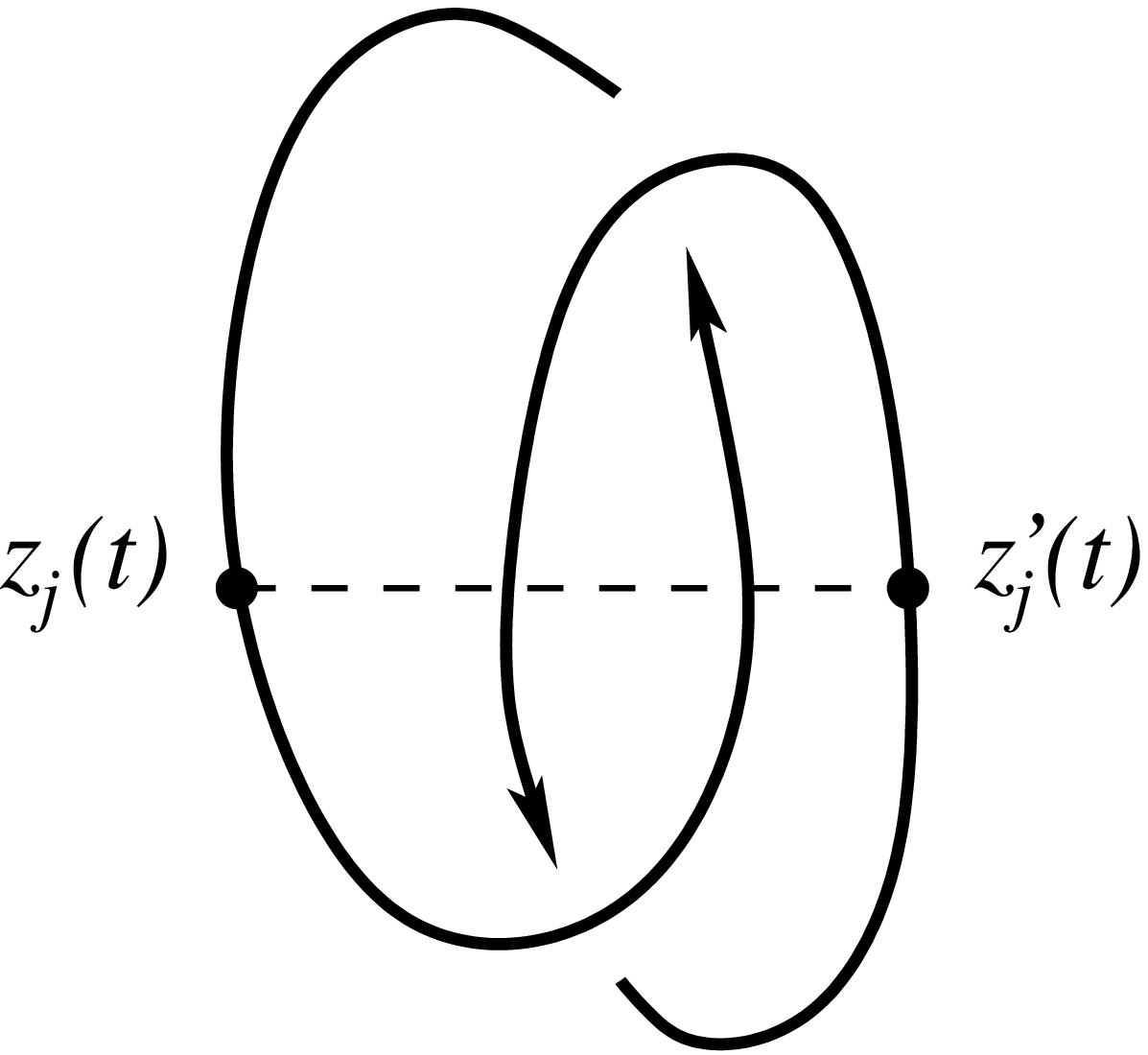}}
\smallskip

The Kontsevich integral can be regarded as a far-going generalization of
these formulas. It aims at encoding all information about how the horizontal
chords on the knot (or tangle) rotate when moved in the vertical direction.
From a more general viewpoint, the Kontsevich integral represents the
monodromy of the Knizhnik--Zamolodchikov connection in the complement
to the union of diagonals in $\C^n$ (see \cite{BN1,Oht}).

\section{Chord diagrams and weight systems}
\label{ch_diag}

\subsection{Algebras $\A(p)$}
The Kontsevich integral of a tangle $T$ takes values in the space of chord
diagrams supported on $T$.

Let $X$ be an oriented one-dimensional manifold, that is, a collection
of $p$ numbered oriented lines and $q$ numbered oriented circles. 
A chord diagram of order $n$
supported on $X$ is a collection of $n$ pairs of unordered points in $X$,
considered up to an orientation- and component-preserving diffeomorphism.
In the vector space formally generated by all chord diagrams of order $n$
we distinguish the subspace spanned by all \textit{four-term relations}
\begin{center}
\includegraphics[width=8cm]{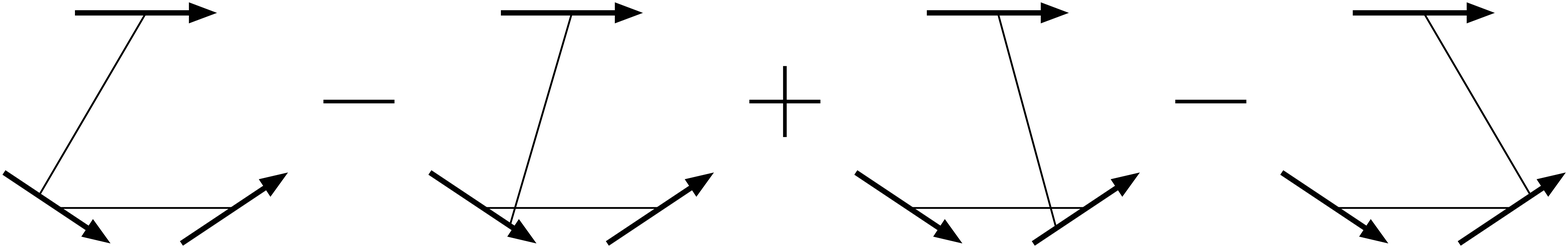}
\end{center}
where thin lines designate chords, while thick lines are pieces
of the manifold $X$. Apart from the shown fragments, all the four diagrams
are identical. The quotient space over all such combinations
will be denoted by $\A_n(X)=\A_n(p,q)$. 
Let $\A(p,q)=\oplus_{n=0}^\infty\A_n(p,q)$
and let $\hat\A(p,q)$ be the graded completion of $\A(p,q)$ (i.e. the space of 
formal infinite series $\sum_{i=0}^\infty a_i$ with $a_i\in\A_i(p,q)$.
If, moreover, we divide $\A(p,q)$ by all \textit{framing independence}
relations (any diagram with an isolated chord, i.e. a chord joining
two adjacent points of the same connected component of $X$, is set to 0),
then the resulting space is denoted by $\A'(p,q)$, and its graded completion 
by $\hat\A'(p,q)$.

The spaces $\A(p,0)=\A(p)$ have the structure of an algebra (the product of
chord diagrams is defined by concatenation of underlying manifolds in
agreement with the orientation). Closing a line component into a
circle, we get a linear map $\A(p,q)\to\A(p-1,q+1)$ which is an isomorphism
when $p=1$. In particular, $\A(S^1)\cong\A(\R^1)$ has the structure of an
algebra; this algebra is denoted simply by $\A$; the Kontsevich integral 
of knots takes its values in its graded completion $\hat\A$. 
Another algebra of special importance is $\hat\A(3)=\hat\A(3,0)$, because it
is where the Drinfeld associators live.

\subsection{Hopf algebra structure}
The algebra $\A(p)$ has a natural structure of a Hopf algebra with
the coproduct $\delta$ defined by all ways to split the set of chords
into two disjoint parts. To give a 
convenient description of its primitive space, one can use
generalized chord diagrams. We now allow 
trivalent vertices not belonging to
the supporting manifold and use STU relations
\begin{center}
\includegraphics[width=5cm]{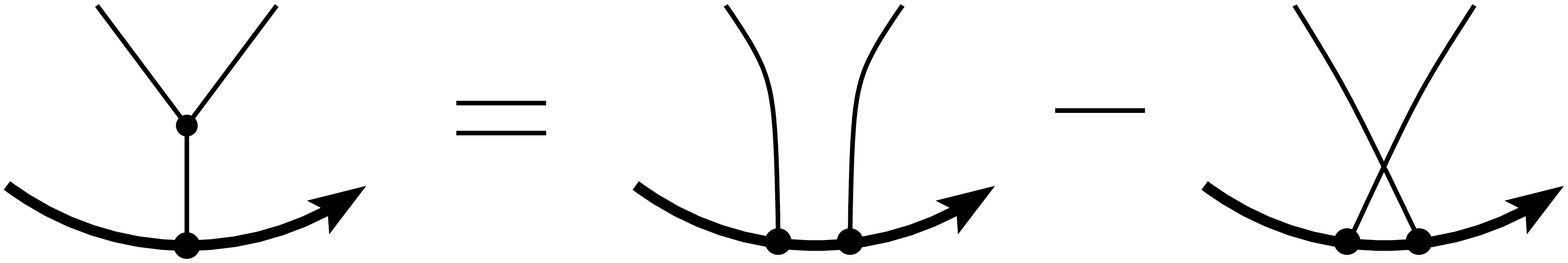}
\end{center}
to express the generalized diagrams as linear combinations of conventional
chord diagrams, e.g.
\begin{center}
\includegraphics[width=6cm]{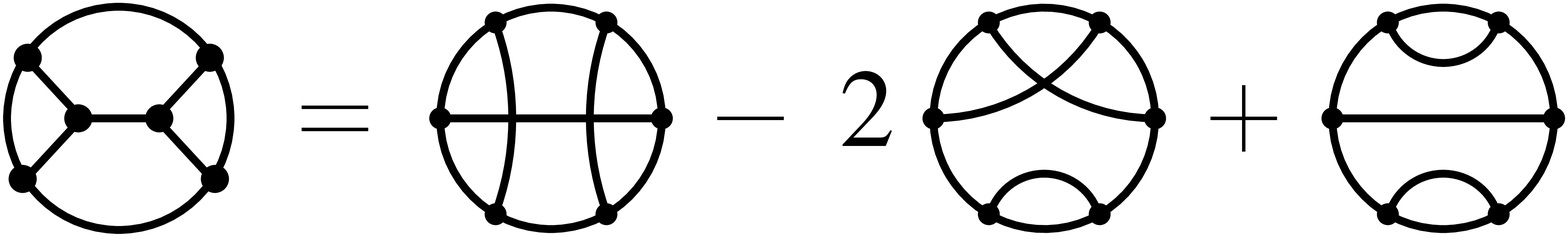}
\end{center}

Then the primitive space coincides with the subspace of $\A(p)$ spanned by
all connected generalized chord diagrams (\textit{connected} means that they
remain connected when the supporting manifold $X$ is disregarded).

\subsection{Weight systems}
A \textit{weight system} of degree $n$ is a linear function on the
space $\A_n$. Every Vassiliev invariant $v$ of degree $n$ defines a weight
system  $\mbox{symb}(v)$ of the same degree called its \textit{symbol}.

\subsection{Algebras $\B(p)$}
Apart from the spaces of chord diagrams modulo four-term relations,
there are closely related spaces of Jacobi diagrams.
A \textit{Jacobi diagram} is defined as a
uni-trivalent graph, possibly disconnected, having at least one vertex of
valency 1 in each connected component and supplied with two additional
structures: a cyclic order of edges
in each trivalent vertex and a labelling of univalent vertices taking values
in the set $\{1,2,\dots,p\}$.
The space $\B(p)$ is defined
as the quotient of the vector space formally generated by all 
$p$-coloured Jacobi diagrams modulo the two types of relations:

\begin{center}
Antisymmetry: \raisebox{-5mm}{\includegraphics[height=11mm]{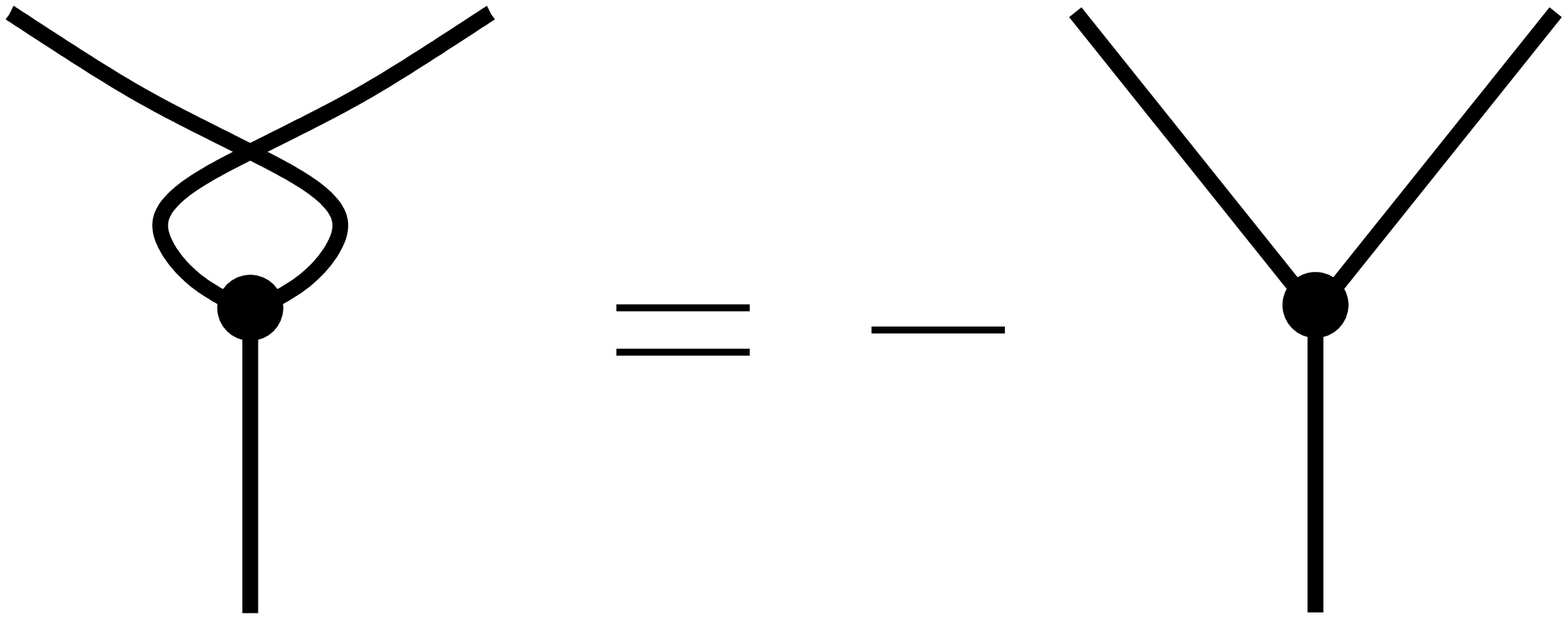}}
\qquad
IHX: \raisebox{-5mm}{\includegraphics[height=11mm]{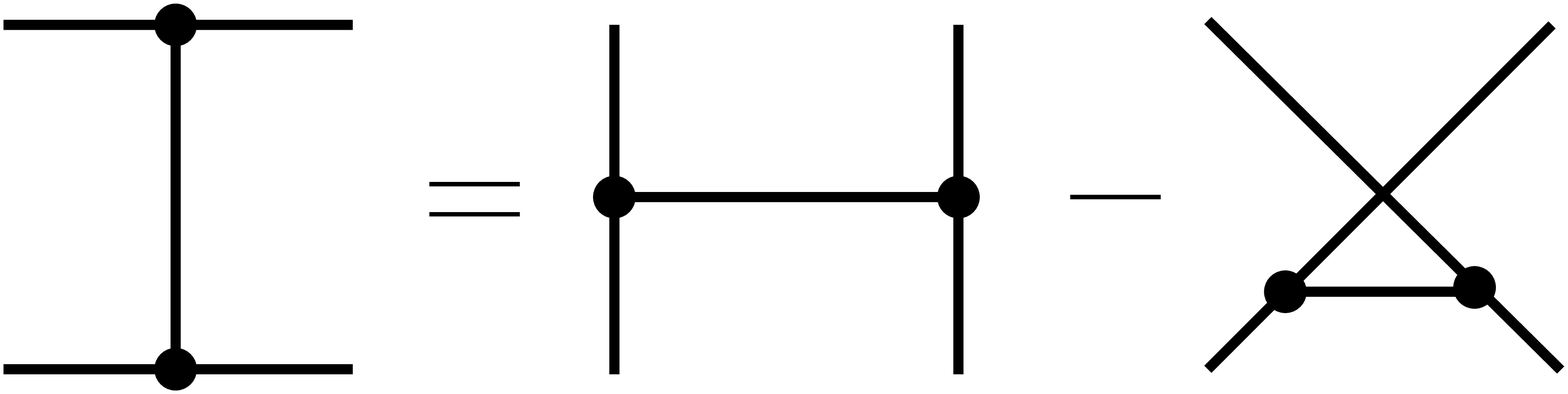}}
\end{center}

The disjoint union of Jacobi diagrams makes the space $\B(p)$ into an
algebra.

The symmetrization map $\chi_p:\B(p)\to\A(p)$, defined as the average
over all ways to attach the legs of colour $i$ to $i$-th connected component
of the underlying manifold:

\begin{center}
\raisebox{-7mm}{\includegraphics[height=15mm]{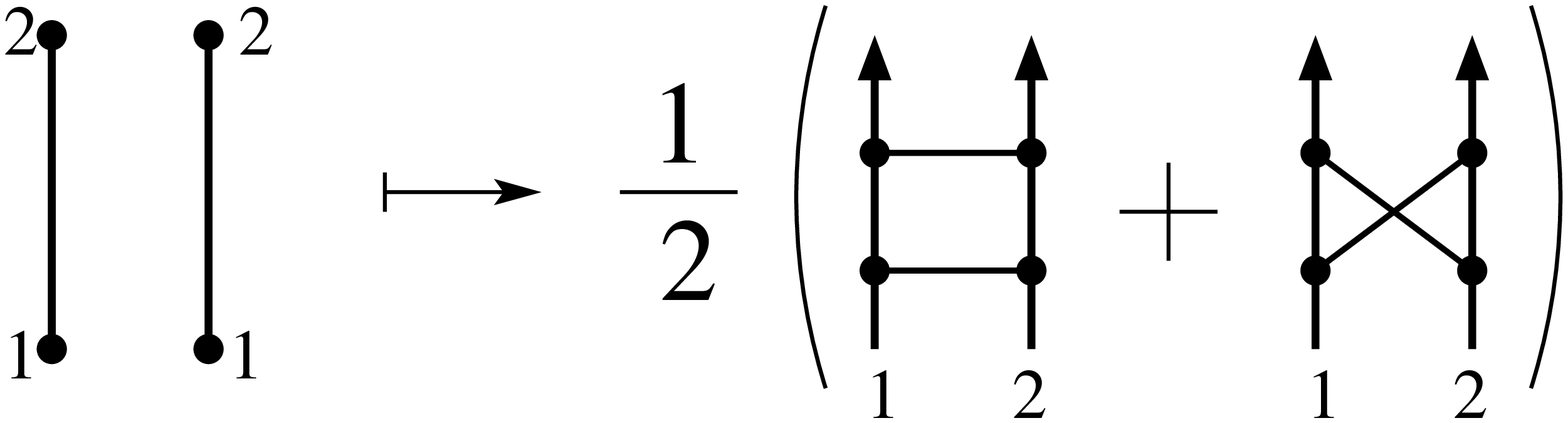}}
\end{center}
is an isomorphism of vector spaces (the formal PBW isomorphism \cite{BN1,LM2}) 
which is not
compatible with the multiplication. The relation between $A(p)$ and $\B(p)$
very much resembles the relation between the universal enveloping algebra
and the symmetric algebra of a Lie algebra.
The algebra $\B=\B(1)$ is used to write out the explicit formula for the
Kontsevich integral of the unknot (see \cite{BLT} and below).

\section{The construction}
\label{defki}

\subsection{Kontsevich's formula}
We will explain the construction of the Kontsevich integral in the classical
case of (closed) oriented knots; for an arbitrary tangle $T$ the formula is 
the same, only the result is interpreted as an element of $\hat\A(T)$.
As above, represent three-dimensional space $\R^3$ as a direct product of a
complex line $\C$ with coordinate $z$ and a real line $\R$ with coordinate
$t$.
The integral is defined for Morse knots,
i.~e. knots $K$ embedded in $\R^3=\C_z\times\R_t$
in such a way that the coordinate $t$ restricted to $K$
has only nondegenerate (quadratic) critical points. (In fact, this condition
can be weakened, but the class of Morse knots is broad enough and 
convenient to work with.)

The  Kontsevich integral $Z(K)$ of the knot $K$ is the following element
of the completed algebra $\hat\A'$:
$$
  Z(K) = \sum_{m=0}^\infty \frac{1}{(2\pi i)^m}
         \int\limits_{\substack{
    t_{\mbox{\tiny min}}<t_m<\dots<t_1<t_{\mbox{\tiny max}}\\
    t_j\mbox{\tiny\ are noncritical}}}\quad
         \sum_{P=\{(z_j,z'_j)\}} (-1)^{\downarrow_P} D_P
         \bigwedge_{j=1}^m \frac{dz_j-dz'_j}{z_j-z'_j} \ .
$$

\subsection{Explanation of the constituents}

The real numbers $t_{\mbox{\scriptsize min}}$ and 
$t_{\mbox{\scriptsize max}}$ are the minimum and the maximum of 
the function $t$ on $K$.

The integration domain is the $m$-dimensional simplex
$t_{\mbox{\scriptsize min}}<t_m<\dots<t_1<t_{\mbox{\scriptsize max}}$ 
divided by the critical values into a certain number of 
{\it connected components}. For example,
the following picture shows an embedding of the unknot where, for $m=2$,
the integration domain has six connected components:
 \begin{center}
     \begin{picture}(280,95)(0,0)
      \put(0,0){\includegraphics[width=280pt]{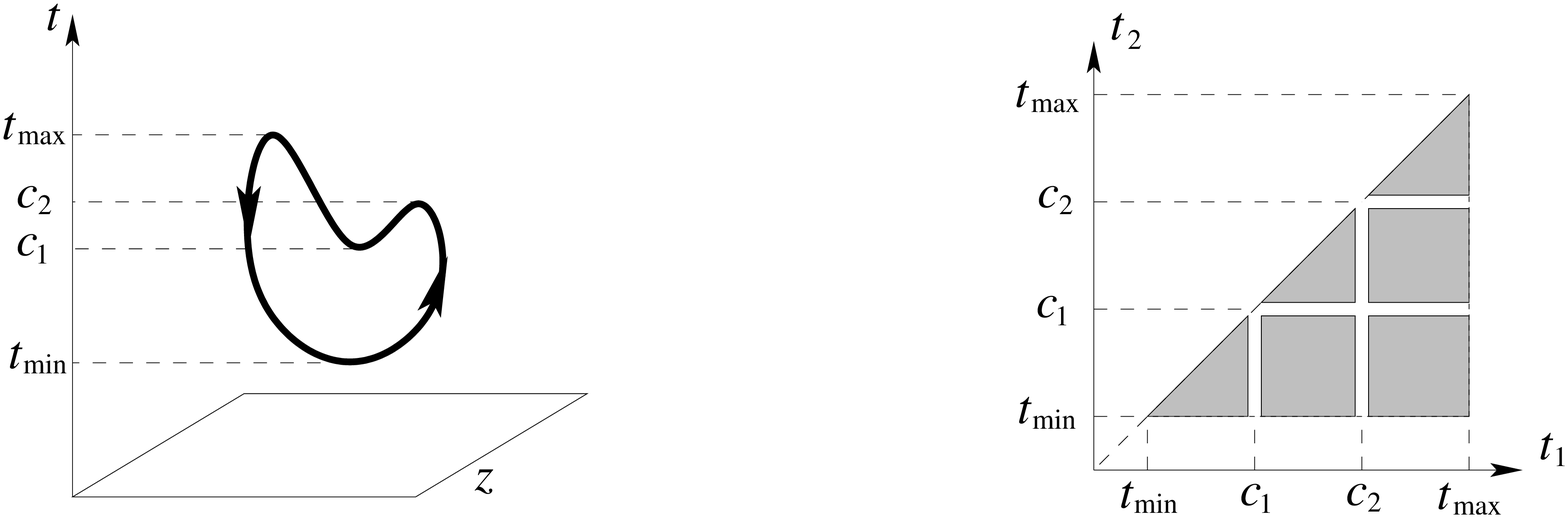}}
     \end{picture}
  \end{center}

The number of summands in the integrand is constant in each connected 
component of the integration domain, but can be different for different 
components.
In each plane $\{t=t_j\}\subset\R^3$ choose an unordered pair of distinct
points $(z_j,t_j)$ and $(z'_j,t_j)$ on $K$, so that $z_j(t_j)$
and $z'_j(t_j)$ are continuous branches of the knot.
We denote by $P=\{(z_j,z'_j)\}$ the collection of such pairs for $j=1,\dots,m$.
The integrand is the sum over all choices of the pairing $P$.
In the example above for the component
$\{t_{\mbox{\scriptsize min}}<t_1<c_1$, 
$c_2<t_2<t_{\mbox{\scriptsize max}}\}$
we have only one possible pair of points on the levels $\{t=t_1\}$ and
$\{t=t_2\}$. Therefore, the sum over $P$ for this component consists of
only one summand. Unlike this, in the component
$\{t_{\mbox{\scriptsize min}}<t_1<c_1$, $c_1<t_2<c_2\}$
we still have only one possibility for the level $\{t=t_1\}$, but the plane
$\{t=t_2\}$ intersects our knot $K$ in four points. So we have
$\binom{4}{2}=6$ possible pairs $(z_2,z'_2)$ and the total number
of summands is six (see the picture below).

For a pairing $P$ the symbol `$\downarrow_P$' denotes the number of points
$(z_j,t_j)$ or $(z'_j,t_j)$ in $P$ where the coordinate $t$ decreases
along the orientation of $K$.

Fix a pairing $P$. Consider the knot $K$ as an oriented circle and connect
the points $(z_j,t_j)$ and $(z'_j,t_j)$ by a chord. Up to a diffeomorphism, 
this chord does not depend on the value of $t_j$ within a connected
component. We obtain a chord diagram with $m$ chords.
The corresponding element of the algebra $\A'$ is denoted by $D_P$.
In the picture below, for each connected component in our example, we show
one of the possible pairings, the corresponding chord diagram with the sign
$(-1)^\downarrow_P$ and the number of summands of the integrand (some of which
are equal to zero in $\A'$ due to the framing independence relation).

\begin{center}
  \includegraphics[width=12cm]{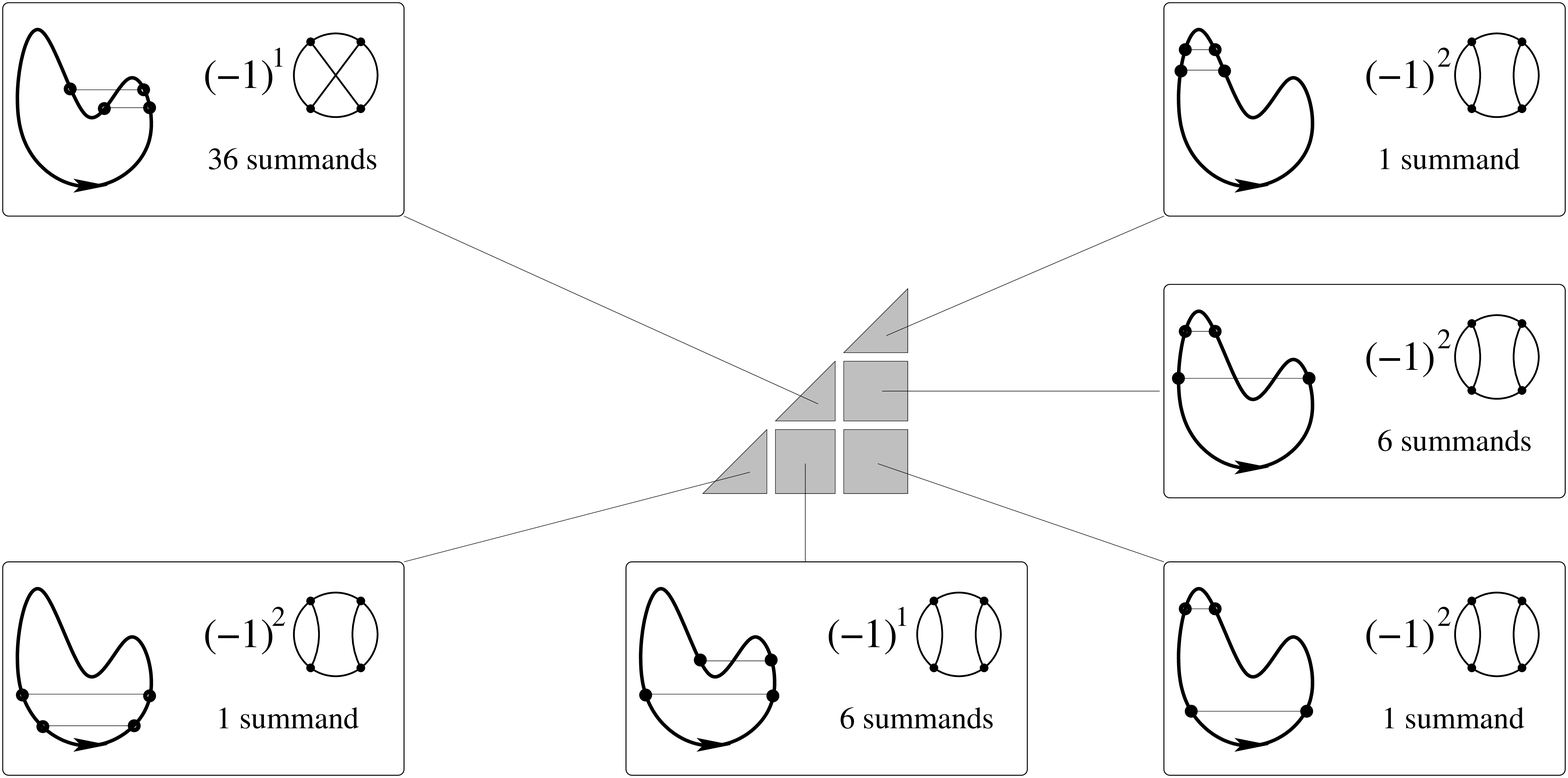}
\end{center}

Over each connected component, $z_j$ and $z'_j$ are smooth functions
of $t_j$.\\ By
${\displaystyle \bigwedge_{j=1}^m \frac{dz_j-dz'_j}{z_j-z'_j}}$ we mean the
{pullback} of this form to the integration domain of variables
$t_1,\dots,t_m$. The integration domain is considered with the 
{orientation} of the space $\R^m$ defined by the natural order of the 
coordinates $t_1$, ..., $t_m$.

By convention, the term in the Kontsevich integral corresponding to $m=0$
is the (only) chord diagram of order 0 with coefficient one. It represents
the unit of the algebra $\A'$.

\subsection{Framed version of the Kontsevich integral}

Let $K$ be a framed oriented Morse knot with writhe number $w(K)$.
Denote the corresponding knot without framing by $\bar{K}$.
The framed version of the Kontsevich integral can be defined
by the formula
$$
 Z^{\mbox{\scriptsize fr}}(K) = 
 e^{\frac{w(K)}{2}\Theta}\cdot Z(\bar{K})\in\hat\A,
$$
where $\Theta$ is the chord diagram with one chord and 
the integral  $Z(\bar{K})\in\hat\A'$ is
understood as an element of the completed algebra
$\hat\A$
 (without 1-term relations) by virtue of a natural
inclusion $\A'\to\A$ defined as identity on the primitive subspace
of $\A'$.  

See \cite{Gor,LM2} for other approaches.

\section{Basic properties}
\label{prop}

\subsection{Constructing the universal Vassiliev invariant}
The Kontsevich integral $Z(K)$
\begin{enumerate}
\item
    converges for any Morse knot $K$,
\item
   is invariant under deformations of the knot in the class 
of Morse knots,
\item
   behaves in a predictable way under the deformation that adds a pair
of new critical points to a Morse knot:
$$
  Z\biggl(\,\raisebox{-12pt}{\includegraphics[width=30pt]{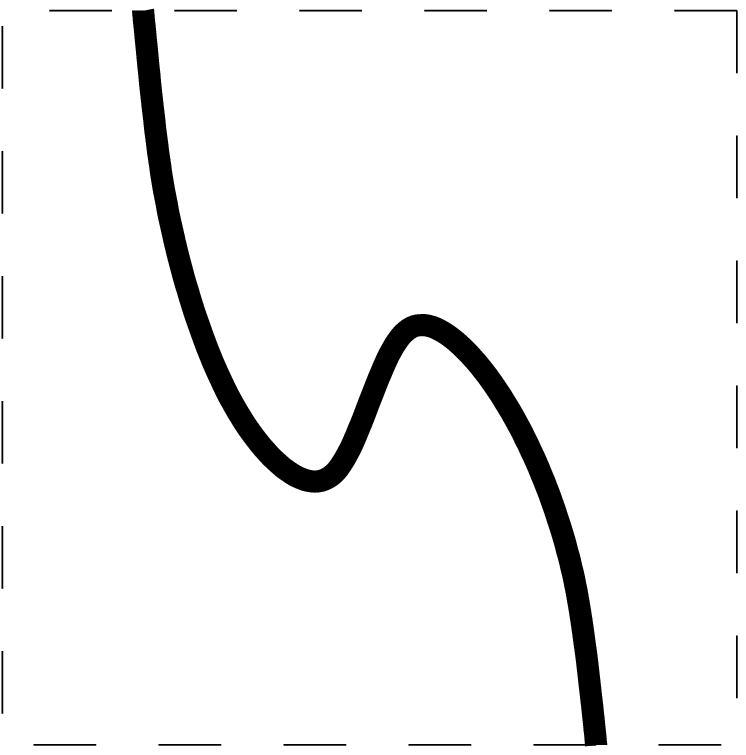}}\,\biggr)
  = Z(H) \cdot 
  Z\biggl(\,\raisebox{-12pt}{\includegraphics[width=30pt]{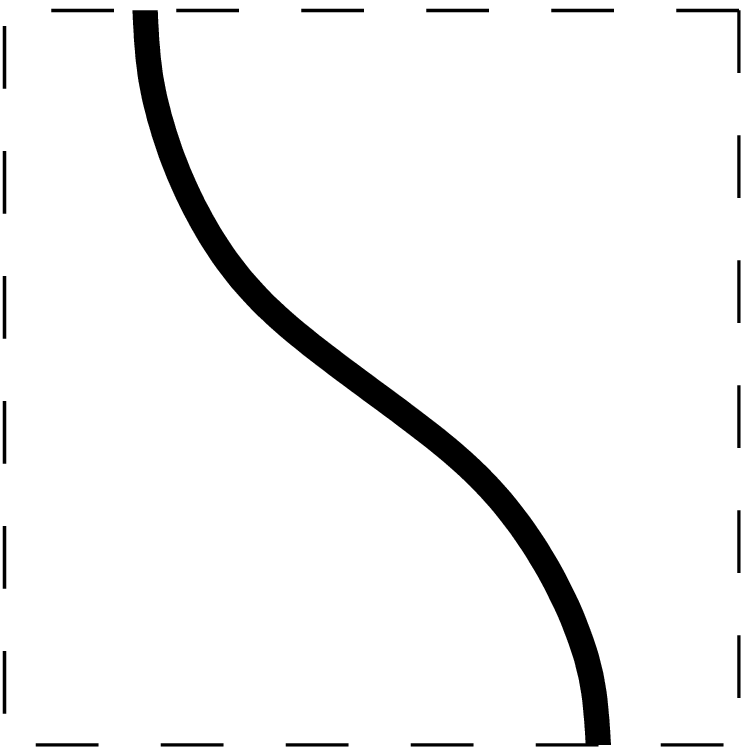}}\,\biggr).
$$
\end{enumerate}

Here the first and the third pictures depict two embeddings of an 
arbitrary knot, differing only in the shown fragment, 
$H = \raisebox{-3mm}{\includegraphics[width=7mm]{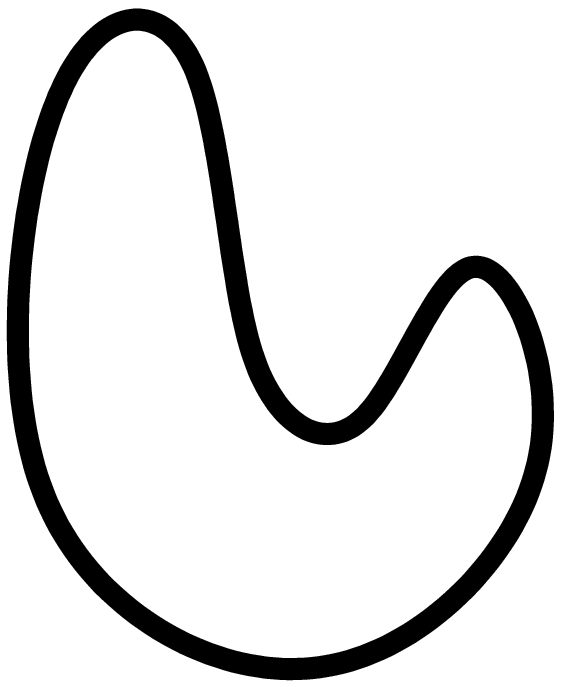}}$
is the {\it hump} (unknot embedded in $\R^3$ in the specified way),
and the product is the product in the completed algebra $\hat\A'$ 
of chord diagrams. The last equality allows one to define a genuine knot
invariant by the formula 
$$\label{I}
  I(K) = Z(K)/Z(H)^{c/2},
$$
where $c$ denotes the number of critical points of $K$ and the ratio
means the division in the algebra $\hat\A'$ according to the rule
$(1+a)^{-1} = 1-a+a^2-a^3+\dots$

The expression $I(K)$ is sometimes referred to as the `final' Kontsevich
integral as opposed to the `preliminary' Kontsevich integral
$Z(K)$. It represents a universal Vassiliev invariant in the following
sense:
\textit{Let $w$ be a weight system, i.e. a linear functional
on the algebra $\hat\A'$. Then the composition $w(I(K))$
is a numerical Vassiliev invariant, and any Vassiliev invariant can 
be obtained in this way.}

The final Kontsevich integral for framed knots is defined in the same
way, using the hump $H$ with zero writhe number.

\subsection{Is universal Vassiliev invariant universal?}
At present, it is not known whether the Kontsevich integral separates knots,
or even if it can tell the orientation of a knot. However, the
corresponding problem is solved, in the affirmative, in the case of
braids and string links (theorem of Kohno--Bar-Natan (\cite{BN2,Koh}).

\subsection{Omitting long chords}
We will state a technical lemma which is highly important
in the study of the Kontsevich integral. It is used in the proof of
the multiplicativity, in the combinatorial construction etc.

Suppose we have a Morse knot $K$ with a distinguished tangle $T$.

\begin{center}
\includegraphics[width=250pt]{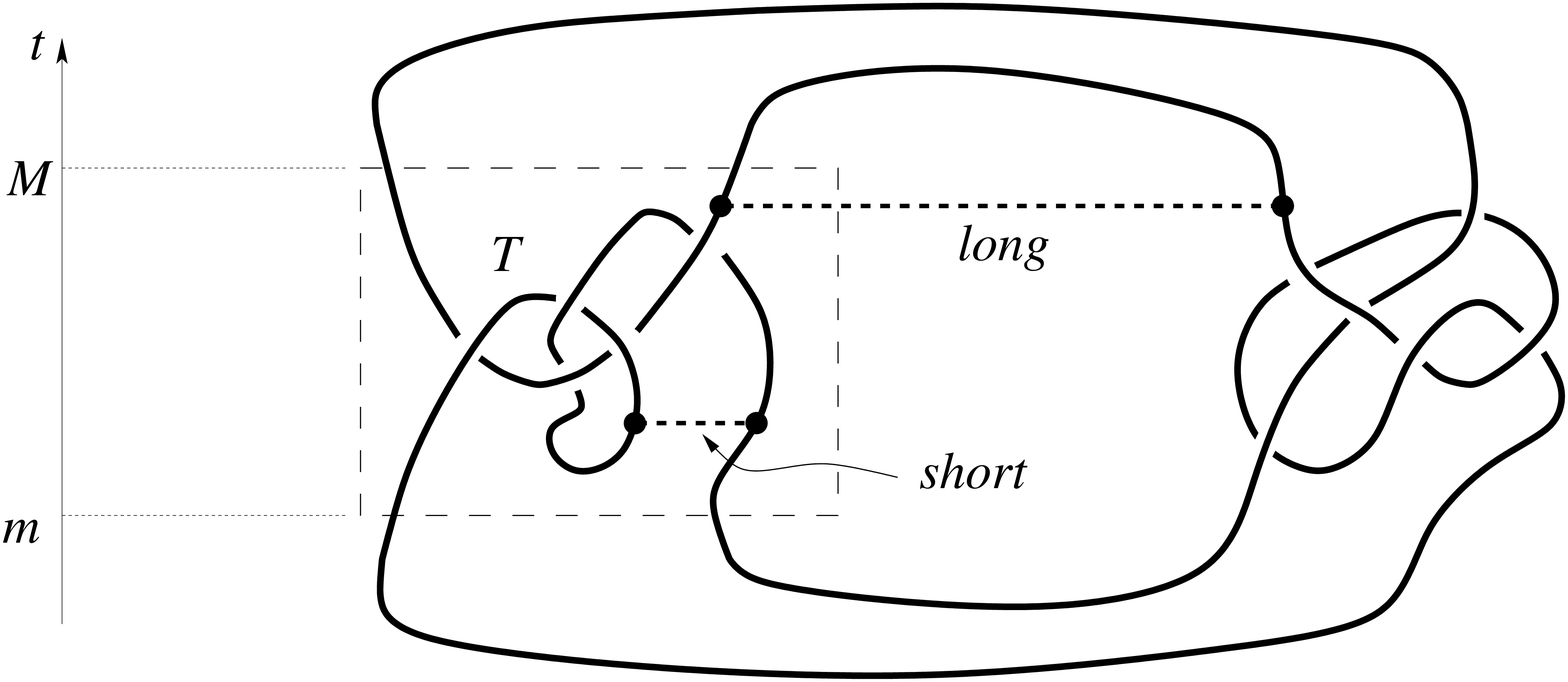}
\end{center}

\noindent
Let $m$ and $M$ be the maximal and
minimal values of $t$ on the tangle $T$. 
In the horizontal planes between the levels $m$ and $M$ we can 
distinguish two kinds of chord: `short' chords that lie either inside
$T$ or inside $K\setminus T$, and `long' chords that connect
a point in $T$ with a point in $K\setminus T$.
Denote by $Z_T(K)$ the expression 
defined by the same formula as the Kontsevich integral $Z(K)$ 
where only short chords are taken into consideration.
More exactly, if $C$ is a
connected component of the integration domain (see Section \ref{defki})
whose projection on the coordinate axis $t_j$ is
entirely contained in the segment $[m,M]$, then in the sum
over the pairings $P$ we
include only those pairings that include short chords.

\textbf{Lemma.}
`Long' chords can be omitted when computing the Kontsevich integral:
$Z_T(K)=Z(K)$.

\subsection{Kontsevich's integral and operations on knots}

The Kontsevich integral behaves in a nice way with respect to the natural
operations on knots, such as mirror reflection,
changing the orientation of the knot, mutation of knots (see \cite{CD}),
cabling (see \cite{Wil1}). We give some details regarding the first two
items.
\smallskip

\textbf{Fact 1.}
Let $R$ be the operation that sends a knot to its mirror image. Define
the corresponding operation $R$ on chord diagrams as
multiplication by $(-1)^n$ where $n$ is the order of the diagram.
Then the Kontsevich integral commutes with
the operation $R$: $Z(R(K))=R(Z(K))$,
where by $R(Z(K))$ we mean simultaneous application of $R$ to all the
chord diagrams participating in $Z(K)$.

\textsl{Corollary.}
The Kontsevich integral $Z(K)$ and the universal Vassiliev invariant 
$I(K)$ of an amphicheiral knot $K$ consist only of even order terms.
(A knot $K$ is called {\it amphicheiral},
if it is equivalent to its mirror image: $K=R(K)$.)
\smallskip

\textbf{Fact 2.}
Let $S$ be the operation on knots which inverts their orientation.
The same letter will also denote the analogous operation on chord
diagrams (inverting the orientation of the outer circle or, which is the 
same thing, axial symmetry of the diagram).
Then the Kontsevich integral commutes with 
the operation $S$ of inverting the orientation:
$Z(S(K))=S(Z(K))$.

\textsl{Corollary.} The following two assertions are equivalent:

--- Vassiliev invariants do not distinguish the orientation of knots,

--- all chord diagrams are symmetric: $D=S(D)$ modulo four-term relations.

The calculations of \cite{Kn} show that up to order 12 all chord 
diagrams are symmetric. For bigger orders the problem is still open.

\subsection{Multiplicative properties}

The Kontsevich integral for tangles is multiplicative:
$$
  Z(T_1)\cdot Z(T_2) = Z(T_1\cdot T_2)
$$
whenever the product $T_1\cdot T_2$, defined by vertical concatenation of
tangles, exists. Here the product in the left-hand side is understood
as the image of the element $Z(T_1)\otimes Z(T_2)$ under the natural map
$\A(T_1)\otimes\A(T_2)\to\A(T_1\cdot T_2)$.

This simple fact has two important corollaries:

\begin{enumerate}
\item
For any knot $K$ the Kontsevich integral $Z(K)$
is a group-like element of the Hopf algebra $\hat\A'$, i.e.
$$
\delta(Z(K)) = Z(K)\otimes Z(K)\,,
$$
where $\delta$ is the comultiplication in $\A$ defined above.
\item
The final Kontsevich integral, taken in a different normalization
$$\label{Iprime}
  I'(K)=Z(H)I(K)=\frac{Z(K)}{Z(H)^{c/2-1}}.
$$
is multiplicative with respect to the connected sum of knots:
$$
  I'(K_1\#K_2)=I'(K_1)I'(K_2),
$$
\end{enumerate}

\subsection{Arithmetical properties}

For any knot $K$ the coefficients in the expansion of $Z(K)$
over an arbitrary basis consisting of chord diagrams are rational 
(see \cite{Kon,LM2} and below).

\section{Combinatorial construction of the Kontsevich integral}
\label{comb}

\subsection{Sliced presentation of knots} 
The idea is to cut the knot into a number of standard simple
tangles, compute the Kontsevich integral for each of them and then
recover the integral of the whole knot from these simple pieces. 


More exactly, we represent the knot by a family of plane diagrams
continuously depending on a parameter $\eps\in(0,\eps_0)$ and cut
by horizontal planes into a number of slices with the following
properties.

\begin{enumerate}
\item
At every boundary level of a slice (dashed lines in the pictures below) 
the distances between various strings are asymptotically
proportional to different whole powers 
of the parameter $\eps$.
\item
Every slice contains exactly one special event and
several strictly vertical strings which
are farther away (at lower powers of $\eps$) from any string
participating in the event than its width. 
\item
There are three types of special events:

\begin{center}
\begin{tabular}{lrcrc}
min/max: & $m=$ & 
\raisebox{-3mm}{\includegraphics[width=25mm]{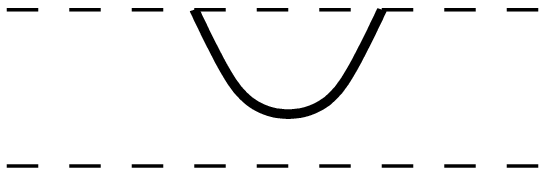}} & $M=$ 
\raisebox{-3mm}{\includegraphics[width=25mm]{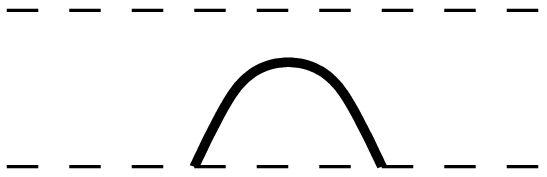}} \\[2em]
braiding:  & $B_+=$ &
\raisebox{-7mm}{\includegraphics[height=15mm]{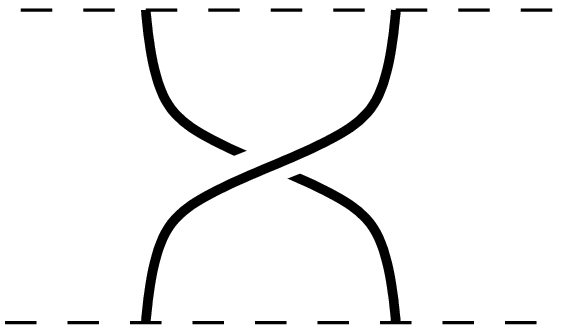}} &  $B_-=$
\raisebox{-7mm}{\includegraphics[height=15mm]{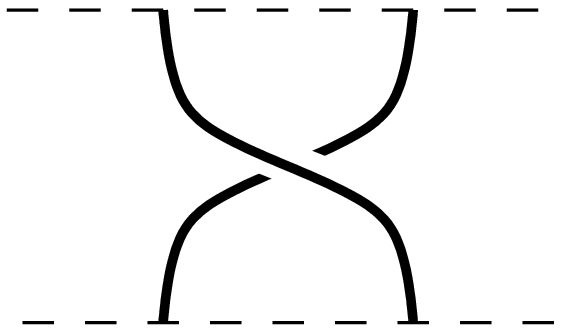}}\\[2em]
associativity: & $A_+=$ &
\raisebox{-7mm}{\includegraphics[height=15mm]{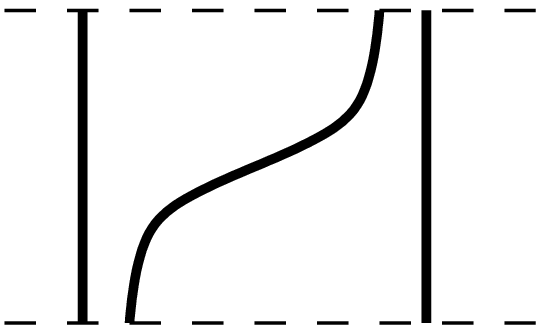}} &  $A_-=$
\raisebox{-7mm}{\includegraphics[height=15mm]{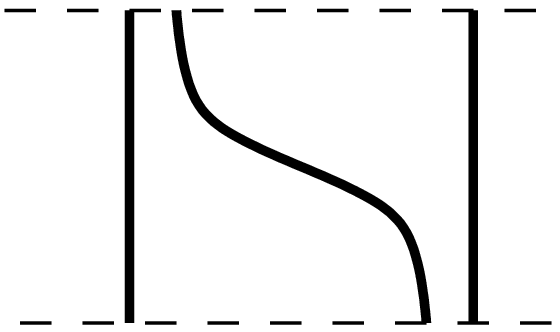}}
\end{tabular}
\end{center}
where, in the two last cases, the strings may be replaced by bunches of
parallel strings which are closer to each other than the width of this
event.
\end{enumerate}


\subsection{Recipe of computation of the Kontsevich integral}

Given such a sliced representation of a knot,
the combinatorial algorithm to compute its Kontsevich integral
consists in the following:
\begin{enumerate}
\item
   Replace each special event by a series of chord diagrams
supported on the corresponding tangle according to the rule
\begin{eqnarray*}
   m,M &\mapsto& 1,\\
   B_+ &\mapsto& R,\quad B_+ \mapsto R^{-1},\\
   A_+ &\mapsto& \Phi,\quad A_- \mapsto \Phi^{-1},
\end{eqnarray*}
where
\begin{align*}
R
  &= \raisebox{-4mm}{\includegraphics[height=11mm]{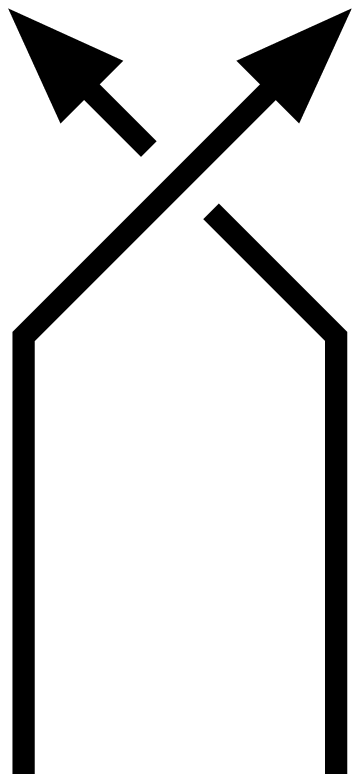}}\cdot
\exp\left(\frac{\raisebox{-1.3mm}{\includegraphics[height=4mm]{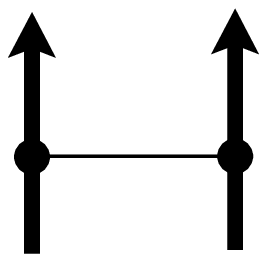}}}{2}\right)
  =  \raisebox{-4mm}{\includegraphics[height=11mm]{R0.eps}}
  + \frac{1}{2}\,\raisebox{-4mm}{\includegraphics[height=11mm]{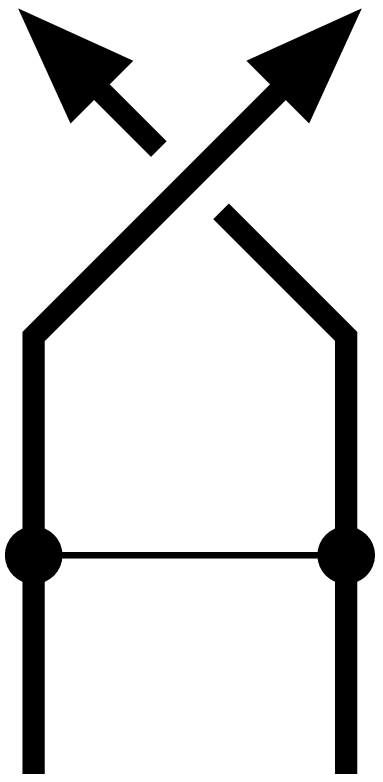}}
  + \frac{1}{2\cdot 2^2}\, \raisebox{-4mm}{\includegraphics[height=11mm]{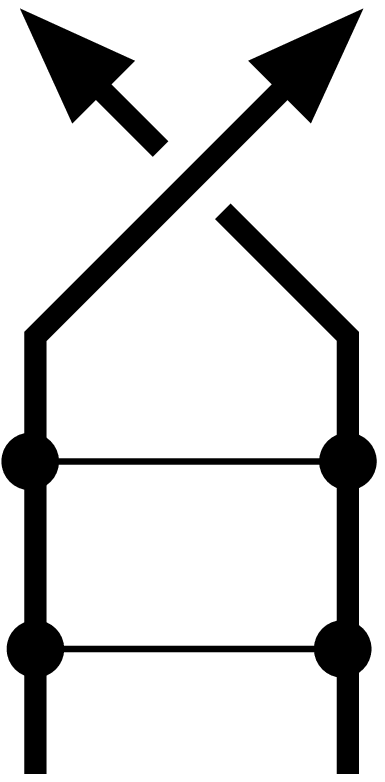}}
  + \frac{1}{3!\cdot 2^3}\, \raisebox{-4mm}{\includegraphics[height=11mm]{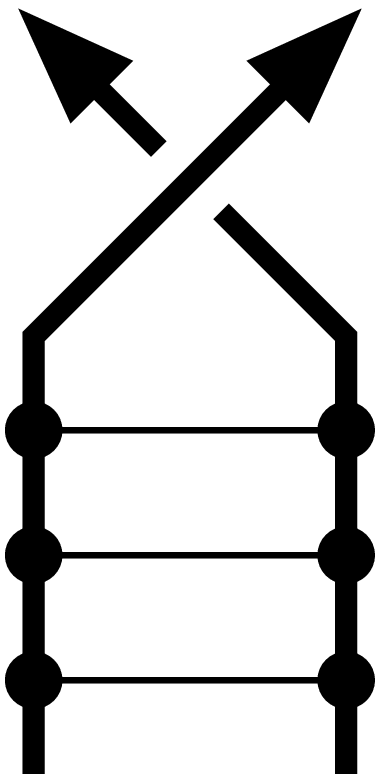}} +\dots\,\\
\Phi &= 1-\frac{\zeta(2)}{(2\pi i)^2} [a,b]
    -\frac{\zeta(3)}{(2\pi i)^3} ([a,[a,b]]+[b,[a,b])+\dots
\end{align*}
($\Phi\in\hat\A(3)$ is the Knizhnik-Zamolodchikov Drinfeld associator 
defined below;
it is an infinite series in two variables
$a=\raisebox{-2mm}{\includegraphics[height=6mm]{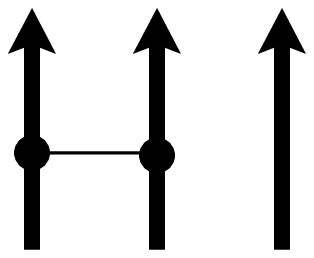}}$, 
$b=\raisebox{-2mm}{\includegraphics[height=6mm]{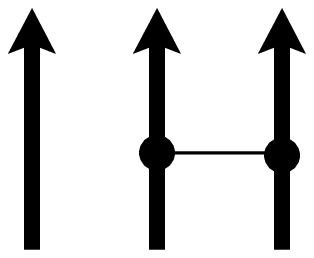}}$).
\item
  Compute the product of all these series from top to bottom taking 
into account the connection of the strands of different tangles, thus 
obtaining an element of the algebra $\hat\A'$.
\end{enumerate}

To accomplish the algorithm, we need two auxiliary operations
on chord diagrams:
\begin{enumerate}
\item
$S_i:\A(p)\to\A(p)$ defined as multiplication by $(-1)^k$
on a chord diagram containing $k$ endpoints of chords on the string
number $i$. This is the correction term in the computation
of $R$ and $\Phi$ in the case when the tangle contains some strings
oriented downwards (the upwards orientation is considered as positive).
\item
$\Delta_i:\A(p)\to\A(p+1)$
acts on a chord diagram $D$ by doubling the $i$-th string of $D$
and taking the sum over all possible lifts of the endpoints of chords
of $D$ from the $i$-th string to one of the two new strings.
The strings are counted by their bottom points from left to right.
This operation can be used to express the combinatorial Kontsevich
integral of a generalized associativity tangle (with strings replaced
by bunches of strings) in terms of the
combinatorial Kontsevich integral of a simple associativity tangle.
\end{enumerate}

\subsection{Example.}
Using the combinatorial algorithm, we compute the Kontsevich integral 
of the trefoil knot $3_1$ to the terms of degree 2.
A sliced presentation for this knot 
shown in the picture implies that
$Z(3_1)= S_3(\Phi)R^{-3}S_3(\Phi^{-1})$
(here the product from left to right corresponds
to the multiplication of tangles from top to bottom).

\noindent
\begin{minipage}{4in}
Up to degree 2, we have:
$\Phi=1+\frac{1}{24}[a,b]+\dots$,
$R=X(1+\frac{1}{2}a+\frac{1}{8}a^2+\dots)$,
where $X$ means that the two strands in each term of the series
must be crossed over at the top.
The operation $S_3$ changes the orientation of the third strand,
which means that $S_3(a)=a$ and $S_3(b)=-b$.
Therefore,
$S_3(\Phi)=1-\frac{1}{24}[a,b]+\dots$,
$S_3(\Phi^{-1})=1+\frac{1}{24}[a,b]+\dots$,
$R^{-3}=X(1-\frac{3}{2}a+\frac{9}{8}a^2+\dots)$
and
$Z(3_1) = (1-\frac{1}{24}[a,b]+\dots)
X(1-\frac{3}{2}a+\frac{9}{8}a^2+\dots)(1+\frac{1}{24}[a,b]+\dots)
=1-\frac{3}{2}Xa-\frac{1}{24}abX
+\frac{1}{24}baX+\frac{1}{24}Xab-\frac{1}{24}Xba+\frac{9}{8}Xa^2+\dots$
Closing these diagrams into the circle, we see that in the algebra $\A$ 
we have $Xa=0$ (by the framing independence relation),
then $baX=Xab=0$ (by the same relation, because these diagrams consist
of two parallel chords) and $abX=Xba=Xa^2=\cdWW$.
The result is:
$Z(3_1)=1+\frac{25}{24}\cdWW+\dots$.
The final Kontsevich integral 
of the trefoil (in the multiplicative normalization, see page
\pageref{Iprime}) is thus equal to
$I'(3_1)=Z(3_1)/Z(H)=(1+\frac{25}{24}\cdWW+\dots)/
(1+\frac{1}{24}\cdWW+\dots)=1+\cdWW+\dots$ 
\end{minipage}
\qquad
\raisebox{-24mm}{\includegraphics[height=55mm]{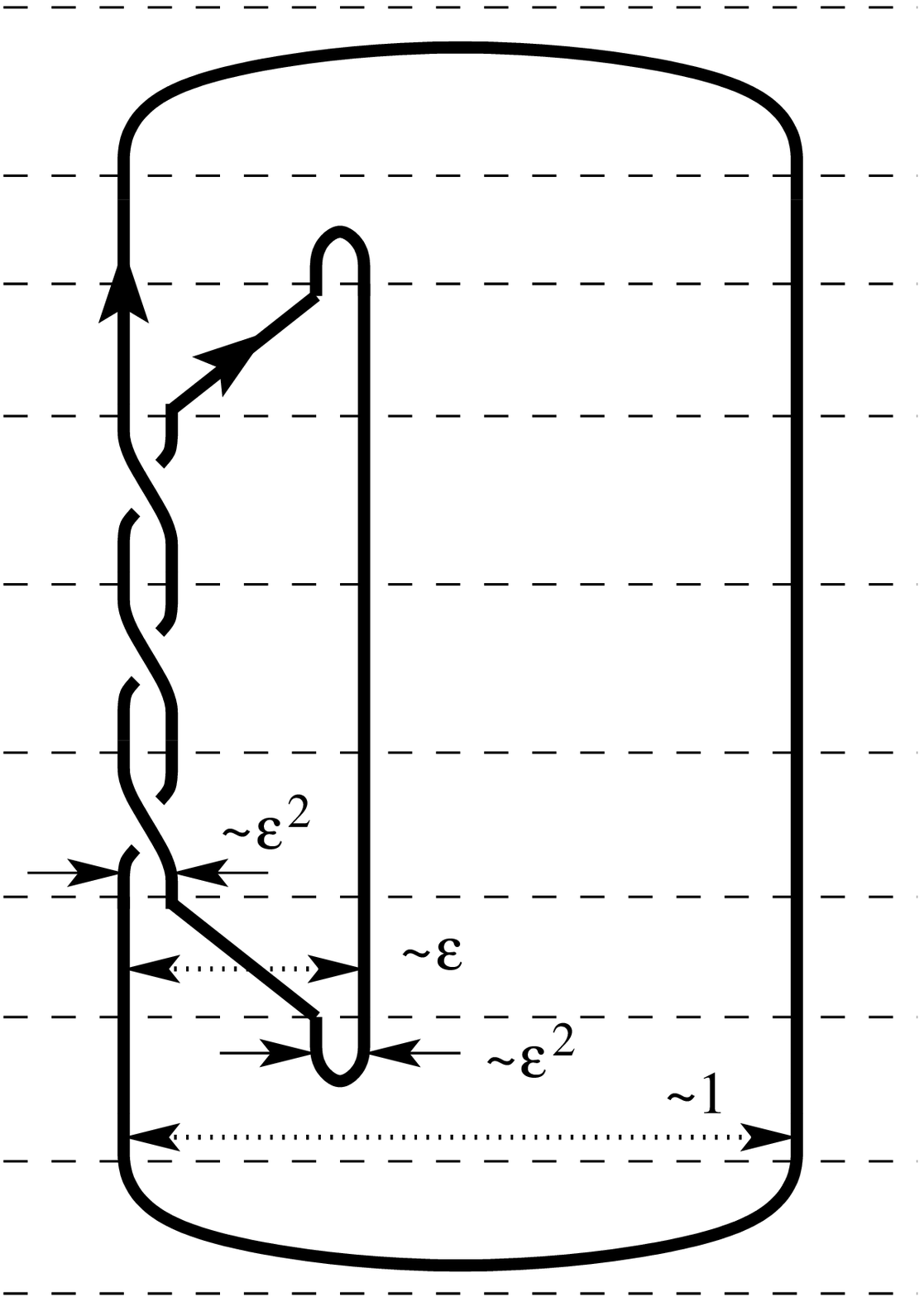}}

\subsection{Drinfeld associator and rationality}

The Drinfeld associator used as a building block in the combinatorial
construction of the Kontsevich integral, can be defined as the limit
$$
\Phi_{\mbox{\scriptsize KZ}}=\lim_{\eps\to0}
\eps^{-b} Z(AT_\eps) \eps^a,
$$
where 
$a=\raisebox{-2mm}{\includegraphics[height=6mm]{vstod.eps}}$,
$b=\raisebox{-2mm}{\includegraphics[height=6mm]{vstdt.eps}}$),
and 
$AT_\eps$ is the positive associativity tangle (special event $A_+$
shown above)
with the distance between the vertical strands constant 1
and the distance between the close endpoints equal to $\eps$.
An explicit formula for $\Phi_{\mbox{\scriptsize KZ}}$ was found by Le and Murakami 
\cite{LM2};
it is written as a nested summation over four variable multiindices and
therefore does not provide an immediate insight into the structure of the
whole series;
we confine ourselves by quoting
the beginning of the series (note that $\Phi_{\mbox{\scriptsize KZ}}$ is a group-like 
element
in the free associative algebra with 2 generators, hence its logarithm
belongs to the corresponding free Lie algebra):
\begin{align*}
  \log(\Phi_{\mbox{\scriptsize KZ}})
 &= -\zeta(2) [x,y]
    -\zeta(3) ([x,[x,y]]+[y,[x,y])\\
 & -\frac{\zeta(2)^2}{10}(4[x,[x,[x,y]]]+[y,[x,[x,y]]]
      +4[y,[y,[x,y]]])\\
 &-\zeta(5) ([x,[x,[x,[x,y]]]]+[y,[y,[y,[x,y]]]])\\
&+(\zeta(2)\zeta(3)-2\zeta(5))([y,[x,[x,[x,y]]]]+ [y,[y,[x,[x,y]]]])\\
&+(\frac{1}{2}\zeta(2)\zeta(3)-\frac{1}{2}\zeta(5)) [[x,y],[x,[x,y]]]\\
&+(\frac{1}{2}\zeta(2)\zeta(3)-\frac{3}{2}\zeta(5)) [[x,y],[y,[x,y]]]+\dots
\end{align*}
where $x=\frac{1}{2\pi i}a$ and $y=\frac{1}{2\pi i}b$.
In general, $\Phi_{\mbox{\scriptsize KZ}}$ is an infinite series whose coefficients
are \textit{multiple zeta values} (\cite{LM2,Zag}) 
$$
  \zeta(a_1,\dots,a_n)=
\sum_{0<k_1<k_2<\dots<k_n} k_1^{-a_1}\dots k_n^{-a_n}.
$$

There are other equivalent definitions of $\Phi_{\mbox{\scriptsize KZ}}$, in particular 
one in terms of the asymptotical behaviour of solutions of the simplest
Knizhnik--Zamolodchikov equation
$$
\frac{dG}{dz} =\left(\frac{a}{z} + \frac{b}{z-1}\right)G,
$$
where $G$ is a function of a complex variable taking values in the
algebra of series in two non-commuting variables $a$ and $b$ 
(see \cite{Dri}).

It turns out (theorem of Le and Murakami \cite{LM2}) that the combinatorial
Kontsevich integral does not change if $\Phi_{\mbox{\scriptsize KZ}}$ is replaced by another
series in $\hat\A(3)$ provided it satisfies certain axioms (among which the
pentagon and hexagon relations are the most important, see \cite{Dri,LM2}).

Drinfeld \cite{Dri} proved the existence of an associator
$\Phi_{\Q}$ with rational coefficients. Using it instead of $\Phi_{\mbox{\scriptsize KZ}}$
in the combinatorial construction, we obtain the following

\textbf{Theorem.} (\cite{LM2}) The coefficients of the Kontsevich integral
of any knot (tangle) are rational when $Z(K)$ is expanded over an arbitrary
basis consisting of chord diagrams.

\section{Explicit formulas for the Kontsevich integral}
\label{explic}

\subsection{The wheels formula}
Let $O$ be the unknot; the expression $I(O)=Z(H)^{-1}$ is referred to
as the \textit{Kontsevich integral of the unknot}. A closed form
formula for $I(O)$ was proved in \cite{BLT}:

\textbf{Theorem.}
$$  I(O) = \exp \sum_{n=1}^\infty b_{2n} w_{2n}
       = 1 + (\sum_{n=1}^\infty b_{2n}w_{2n})
           + \frac{1}{2}(\sum_{n=1}^\infty b_{2n}w_{2n})^2 + \dots
$$
Here $b_{2n}$ are modified Bernoulli numbers, i.e.
the coefficients of the Taylor series
$$
  \sum_{n=1}^\infty b_{2n}x^{2n} = \frac{1}{2}
           \ln\frac{e^{x/2}-e^{-x/2}}{x},
$$
($b_2=1/48$, $b_4=-1/5760$, $b_6=1/362880$,\dots),
and $w_{2n}$ are the {\it wheels}, i.~e. Jacobi diagrams
 of the form
$$
  w_2=\raisebox{-3.5mm}{\includegraphics[height=8mm]{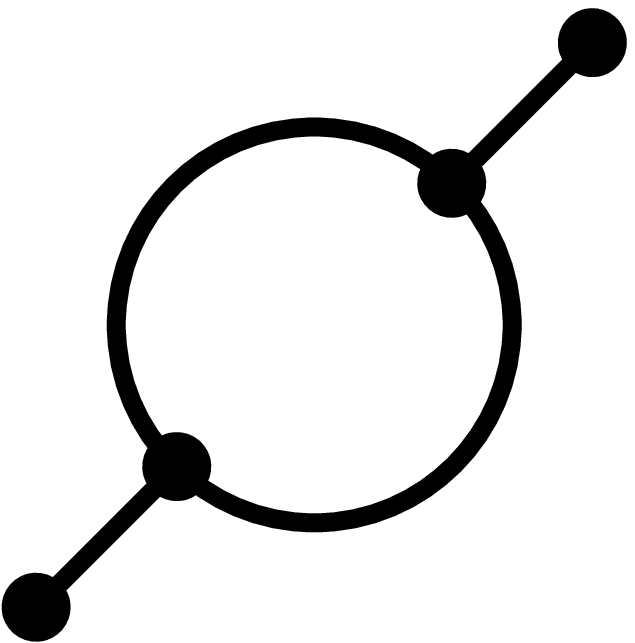}},\quad
  w_4=\raisebox{-3.5mm}{\includegraphics[height=8mm]{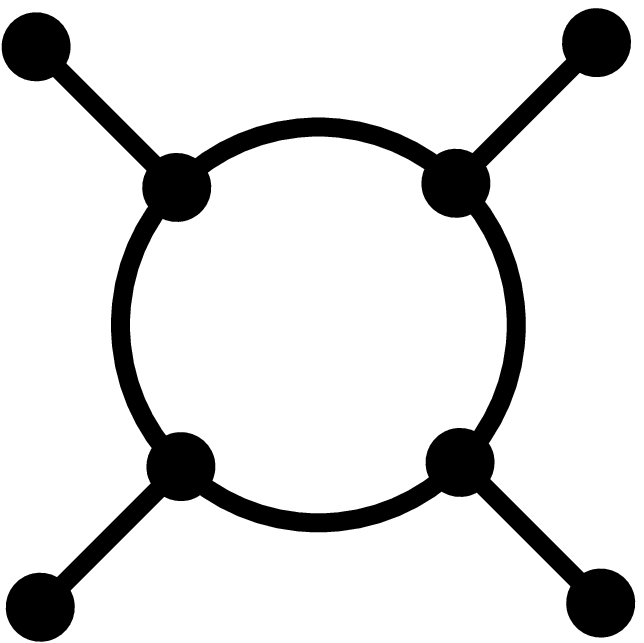}},\quad
  w_6=\raisebox{-5mm}{\includegraphics[height=11mm]{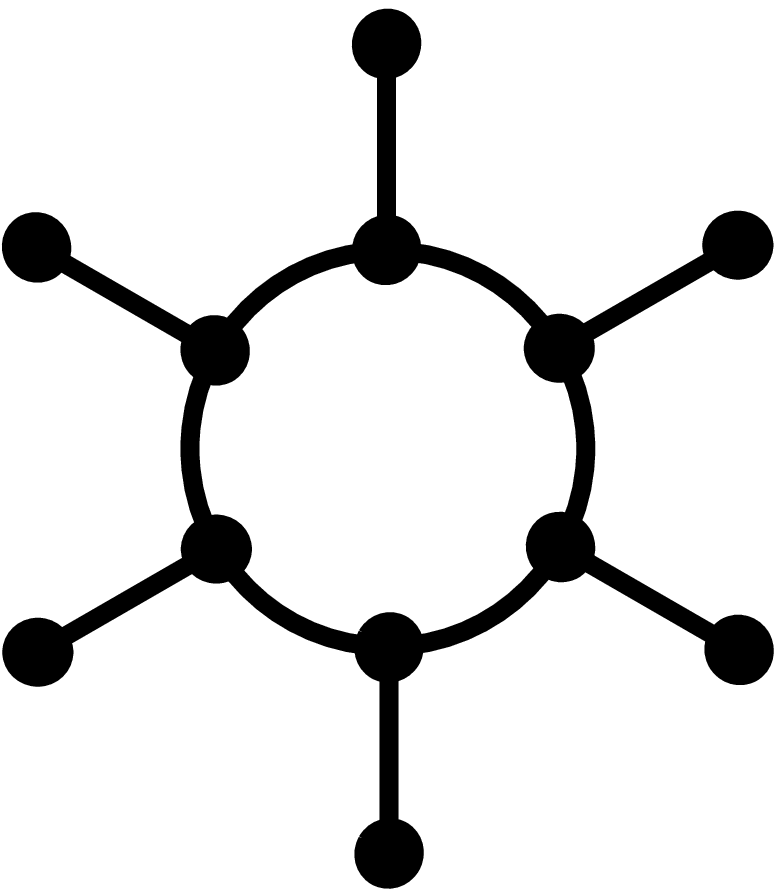}},\quad\dots
$$
The sums and products are understood as operations in the 
algebra of Jacobi diagrams $\B$, and the result is then carried over
to the algebra of chord diagrams $\A$ along the isomorphism $\chi$ 
(see Section \ref{ch_diag}).

\subsection{Generalizations}

There are several generalizations of the wheels formula:

1. Rozansky's \textit{rationality conjecture} \cite{Roz}
proved by A.Kricker \cite{Kri} affirms that the Kontsevich integral
of any (framed) knot can be written in a form resembling the wheels formula.
Let us call the \textit{skeleton} of a Jacobi diagram the regular 3-valent
graph obtained by `shaving off' all univalent vertices.
Then the wheels formula says that all diagrams in the expansion of
$I(O)$ have one and the same skeleton (circle), and the generating function
for the coefficients of diagrams with $n$ legs is a certain analytic
function, more or less rational in $e^x$. In the same way, the theorem of
Rozansky and Kricker states that the terms in $I(K)\in\hat\B$, when arranged
by their skeleta, have the generating functions of the form
$p(e^x)/A_K(e^x)$, where $A_K$ is the Alexander polynomial of $K$ and $p$ is
some polynomial function.  Although this theorem does not give an explicit
formula for $I(K)$, it provides a lot of information about the structure of
this series.

2. J.March\'e \cite{Mar} gives a closed form formula for the
Kontsevich integral of torus knots $T(p,q)$.

The formula of March\'e, although explicit, is rather intricate, and 
here, by way of example, we only write out the first several terms of the 
final Kontsevich integral $I'$ 
for the trefoil (torus knot of type (2,3)), following \cite{Wil2}:
$$
  I'(\raisebox{-3mm}{\includegraphics[height=9mm]{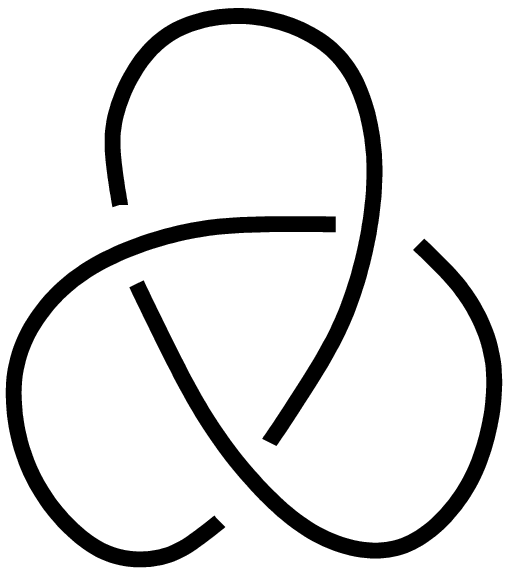}})=
 \raisebox{-2.5mm}{\includegraphics[height=7mm]{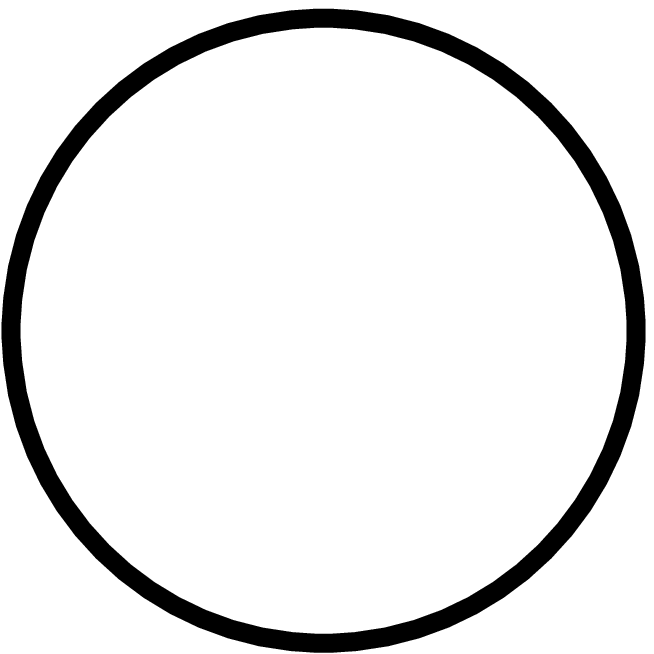}}
-\raisebox{-2.5mm}{\includegraphics[height=7mm]{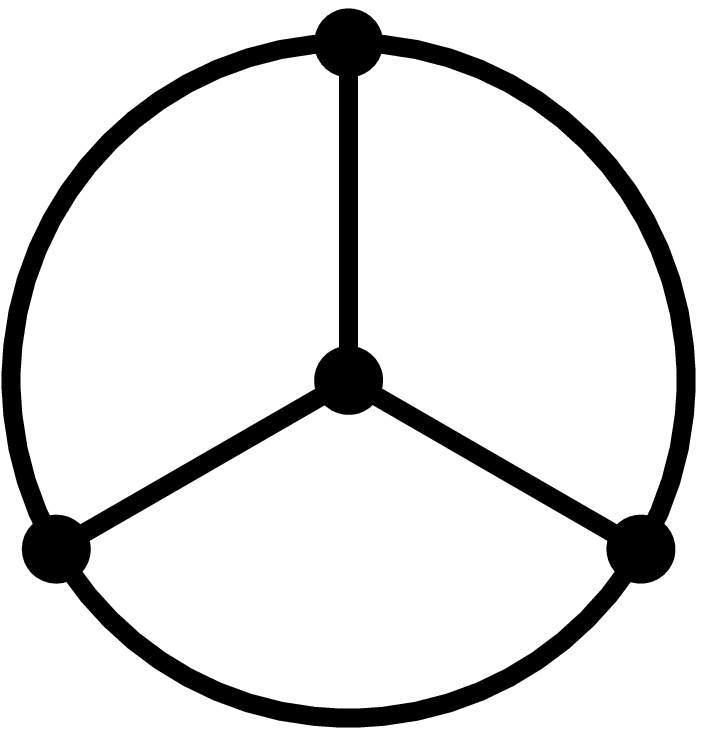}}
+\raisebox{-2.5mm}{\includegraphics[height=7mm]{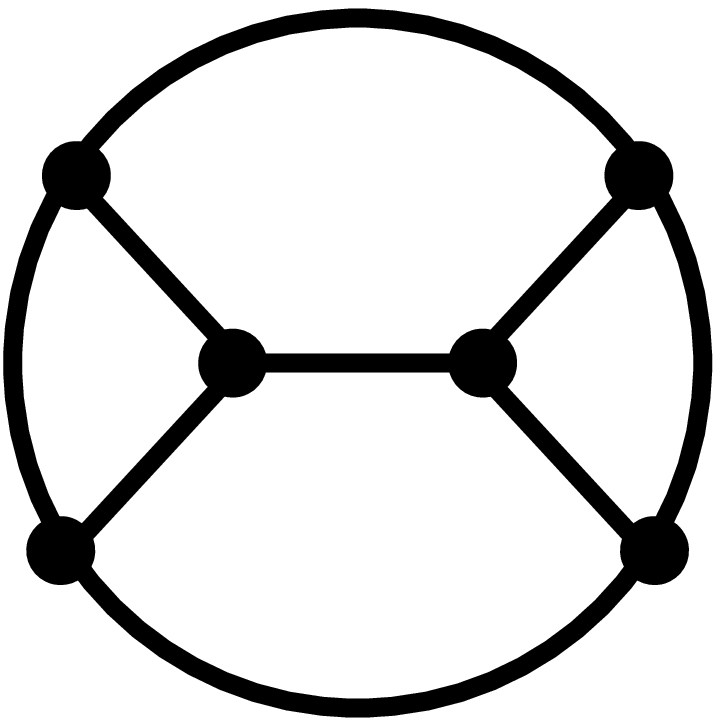}}
-\frac{31}{24}\,\raisebox{-2.5mm}{\includegraphics[height=7mm]{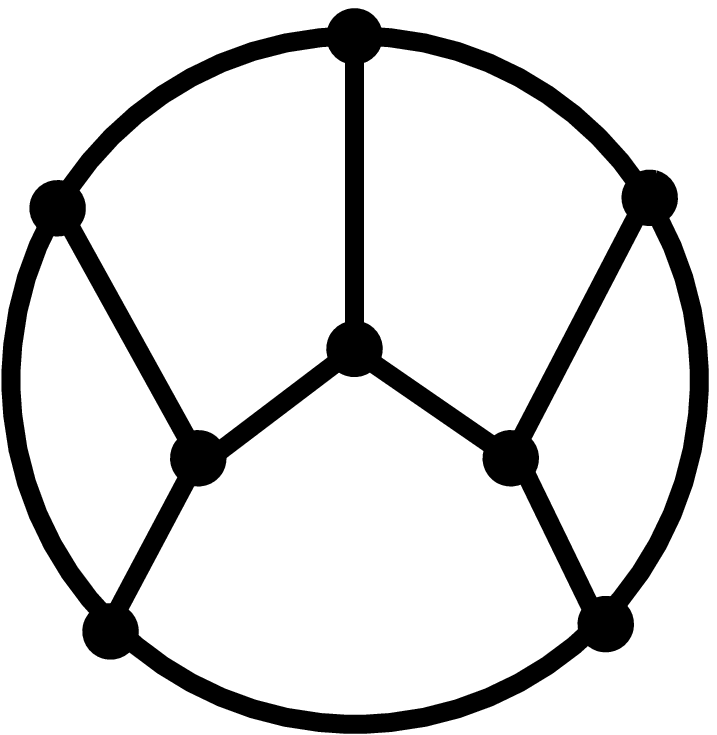}}
+\frac{5}{24}\,\raisebox{-2.5mm}{\includegraphics[height=7mm]{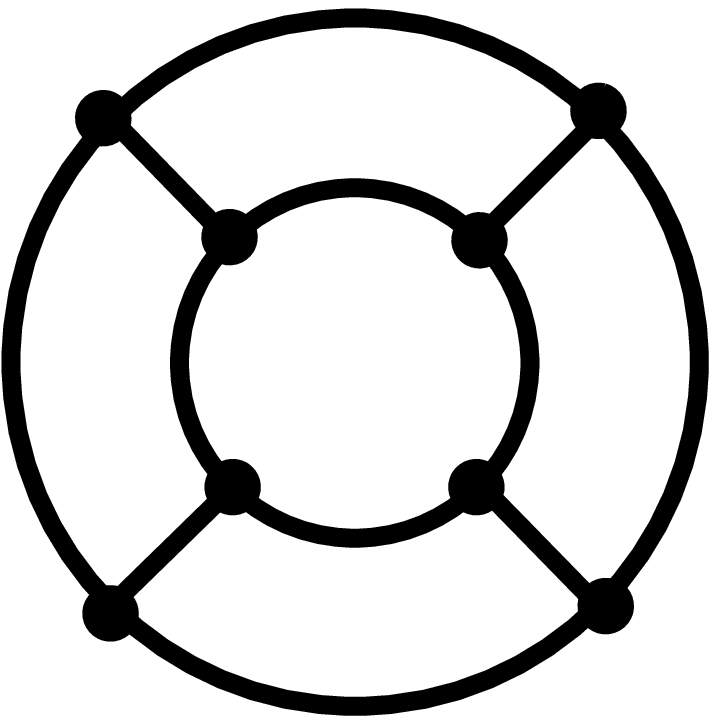}}
+\frac{1}{2}\,\raisebox{-2.5mm}{\includegraphics[height=7mm]{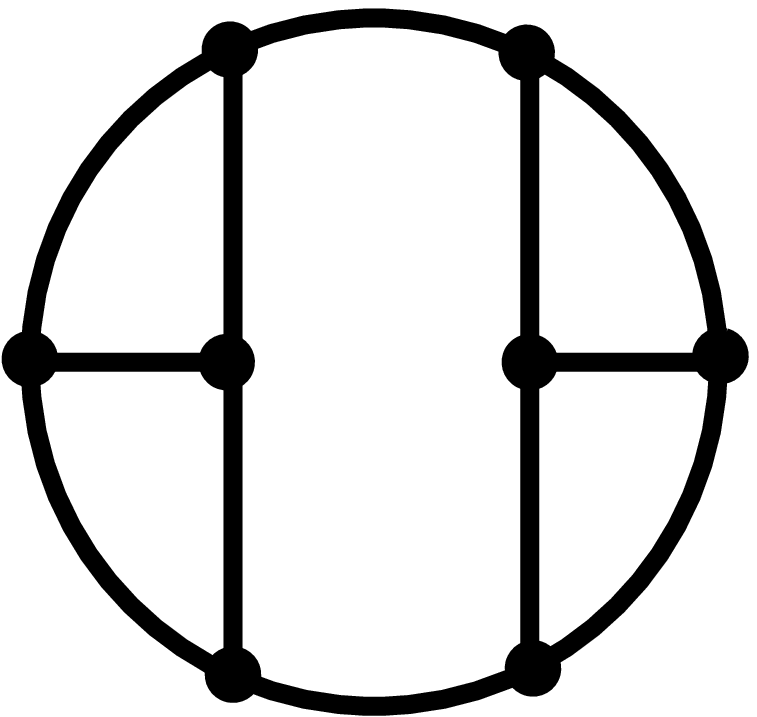}}
+\dots
$$

\subsection{First terms of the Kontsevich integral}

A Vassiliev invariant $v$ of degree $n$ is called 
\textit{canonical} if it can be recovered from the Kontsevich integral
by applying a homogeneous weight system, i.e. if
$v = \mbox{symb}(v)\circ I$.
Canonical invariants define a grading in the filtered space of Vassiliev
invariants which is consistent with the filtration.
If the Kontsevich integral is expanded over a fixed basis
in the space of chord diagrams $\hat\A'$, then the coefficient
of every diagram is a canonical invariant.
According to \cite{Sta,Wil2}, the expansion of the final Kontsevich integral
up to degree 4 can be written as follows:
\begin{align*}
  I'(K) &= 
 \raisebox{-2.5mm}{\includegraphics[height=7mm]{cd0.eps}}
-c_2(K)\,\raisebox{-2.5mm}{\includegraphics[height=7mm]{fd2.eps}}
-\frac{1}{6}j_3(K)\,\raisebox{-2.5mm}{\includegraphics[height=7mm]{fd3.eps}}\\
&+\frac{1}{48}\bigl(4j_4(K)+36c_4(K)-36c_2^2(K)+3c_2(K)\bigr)
\,\raisebox{-2.5mm}{\includegraphics[height=7mm]{fd41.eps}}\\
&+\frac{1}{24}\bigl(-12c_4(K)+6c_2^2(K)-c_2(K)\bigr)
\,\raisebox{-2.5mm}{\includegraphics[height=7mm]{fd42.eps}}\\
&+\frac{1}{2}c_2^2(K)\,
\,\raisebox{-2.5mm}{\includegraphics[height=7mm]{fd43.eps}}
+\dots
\end{align*}
where $c_n$ are coefficients of the Conway polynomial 
$\nabla_K(t)=\sum c_n(K)t^n$
and $j_n$ are modified coefficients of the 
Jones polynomial  $J_K(e^t)=\sum j_n(K)t^n$.
Therefore, up to degree 4, the basic canonical Vassiliev invariants
of unframed knots are $c_2$, $j_3$, $j_4$, $c_4+\frac{1}{12}c_2$
and $c_2^2$.

\end{document}